# COMPARISON OF SYMBOLIC AND ORDINARY POWERS OF PARITY BINOMIAL EDGE IDEALS

Nadia Taghipour, Shamila Bayati, Farhad Rahmati

**ABSTRACT.** In this paper, we investigate when symbolic and ordinary powers of the parity binomial edge ideal of a graph fail to be equal. It turns out that if $\mathcal{I}_G$ is the parity binomial edge ideal of a graph $G$, then in each of the following cases the symbolic power $\mathcal{I}_G^{(t)}$ and the ordinary power $\mathcal{I}_G^t$ are not equal for some $t$: (i) the clique number of $G$ is greater than 3; (ii) $G$ has a net; or (iii) $G$ has a PT as an induced subgraph.

**2020 Mathematics Subject Classification**: 13C13, 13F20, 05E40.

**Keywords**: Clique number; Net-free graph; Parity binomial edge ideal; Symbolic power;

## 1 Introduction

Binomial ideals with a wide range of applications found a new way to combinatorial commutative algebra via binomial edge ideal of graphs, as first defined in [19] and [31]. Later, other binomial ideals associated with a graph have been introduced, namely permanental edge ideals in [22] and parity binomial edge ideals in [26]. In this paper, we deal with the latter which are naturally defined associated with a graph as follows: let $G$ be a simple graph on the vertex set $V(G) = [n]$ with the edge set $E(G)$, and suppose that $S = K[x_1, \ldots, x_n, y_1, \ldots, y_n]$ is the polynomial ring over a field $K$ in $2n$ variables. The parity binomial edge ideal of $G$ is the ideal

$$\mathcal{I}_G = (x_i x_j - y_i y_j : \{i, j\} \in E(G)) \subseteq S.$$

There exists a rich literature on the symbolic powers of ideals which reflects their geometric significance. Recall that if $\mathfrak{p}_1, \ldots, \mathfrak{p}_r$ are the minimal primes of an ideal $I$ in a Noetherian ring, then its $m$-th symbolic power is defined to be $I^{(m)} = \mathfrak{q}_1 \cap \cdots \cap \mathfrak{q}_r$ where $\mathfrak{q}_i$ is the $\mathfrak{p}_i$-primary component of $I^m$. While ordinary powers of an ideal $I$ naturally come up in the algebraic context, from the geometric point of view they are the ideals with the same variety as $I$ but not necessarily with the same associated primes as $I$. This defect, from the viewpoint of geometry, makes symbolic powers more important for a geometer; see also [9, Theorem 3.14] and [33, Corollary 2.9] in the frame of geometry for a nice interpretation of the symbolic power $I^{(m)}$ as the set of all polynomials vanishing to order $m$ over the variety of $I$ where $I$ is a homogeneous radical ideal in a polynomial ring over an algebraically closed field $K$ of characteristic zero. On the other hand, even if we have generators of $I$, it is not easy to find a finite set of generators of $I^{(m)}$; an advantage that attracts attention to ordinary powers and makes the comparison between the symbolic and ordinary powers important; in [35], one can see a survey on the so-called *containment problem* which concerns which ordinary powers include a given symbolic power of $I$. See also an overview of related results in [6].

In this paper, we focus on the comparison of symbolic and ordinary powers of parity binomial edge ideals. While there are various results for monomial ideals associated with the combinatorial structures whose symbolic and ordinary powers are equal, less is known about the equality of these powers of binomial ideals in the combinatorial context; regarding monomial ideals, for example, see [14, 15, 18] for a translation of the combinatorial property of the max-flow-min-cut and Conforti-Cornuéjols conjecture in terms of the equality of symbolic and ordinary powers of squarefree monomial ideals. Besides, one can find in [3, 8, 15, 20, 21] how investigating the coincidence of symbolic and ordinary powers of a monomial ideal through its symbolic Rees algebra is related to the associated combinatorial structure. One also may see in [1, 2, 7, 12, 21, 32, 36] how some classes of squarefree monomial ideals with equal symbolic and ordinary powers are characterized in terms of underlying combinatorial structures.

Recently equality of these powers for binomial edge ideals has attracted more attention. Ene and Herzog show in [11, Theorem 3.3] that if $J$ is the binomial edge ideal of a simple connected graph, and $<$ is a monomial order such that the initial ideal $\text{in}_<(J)$ is a squarefree monomial ideal, then the equality of $\text{in}_<(J)^{(m)} = \text{in}_<(J)^m$ for all $m$ implies that $J^{(m)} = J^m$ for all $m$. By this result, we may reduce the problem about powers of binomial edge ideals of some classes of graphs to squarefree monomial ideals and have the support of rich and growing literature on squarefree monomial ideals whose symbolic and ordinary powers are coincident. Another result regarding the comparison of symbolic and ordinary powers of binomial edge ideals is by Ene et al. which they, along with other results, show in [10, Theorem 4.1] that if $G$ is a block graph whose binomial edge ideal $J$ is Cohen-Macauley, then $J^{(m)} = J^m$ for all $m$ if and only if $G$ is net-free; a result which is partially generalized in [24, Theorem 3.11] for generalized caterpillar graphs whose binomial edge ideal is not necessarily Cohen-Macaulay. It is also shown in [30, Theorem 4.3] that if $J$ is the binomial edge ideal of a complete multipartite graph, then $J^{(m)} = J^m$ for all $m$. To the best of your knowledge, these are the only results regarding the comparison of binomial ideals in the combinatorial context.

Now we are going to study the coincidence of symbolic and ordinary powers of parity binomial edge ideals by the structure of the underlying graph. It turns out in Theorem 2.4 that parity binomial edge ideals with coincident symbolic and ordinary powers have a net-free underlying graph. In Theorem 2.6, we also see that if $\mathcal{I}_G^{(t)} = \mathcal{I}_G^t$ for every $t$, then $G$ does not have an induced so-called PT graph, that is, a cycle of length five and a cycle of length 3 sharing an edge; see Section 2 for more explanation.

By [34, Proposition 3.9] and [27, Remark 3.4], if $\text{char } K \neq 2$, and $G$ is the complete graph $K_r$ for some $r \geq 4$, then there exists a positive integer $n$ such that $\mathcal{I}_G^m \neq \mathcal{I}_G^{(m)}$. In Theorem 2.5, generalizing this fact, we show that if $G$ has a clique of size $r$ for some $r \geq 4$, then there exists such a positive integer $m$ with $\mathcal{I}_G^m \neq \mathcal{I}_G^{(m)}$.

The last section of this paper is devoted to some open questions. As pointed out by Morey in [29], the problem of coincidence of all ordinary and symbolic powers of an ideal in a Noetherian ring reduces to check only for a finite number of powers less than or equal to an invariant of the ideal called the index of stability. See also [28, Theorem 4.8] where it is shown that for a squarefree monomial ideal whose minimal set of generators has $\mu$ elements, if $I^{(m)} = I^m$ for all $m \leq \lceil \frac{\mu}{2} \rceil$, then the equality holds for all $m$. Even more, sometimes it happens that the equality $I^{(m)} = I^m$ for some $m$ is enough to have equality for all $m$. As a result of [23, Corollary 2.5], if $P$ is a prime ideal of height two in a regular local ring $R$ of dimension 3, one has the equality of all ordinary and symbolic powers of $P$ if and only if $P^{(m)} = P^m$ for some $m > 1$; see [17] for more explanation on this result. A similar result for the saturated homogeneous ideal defining a subscheme in $\mathbb{P}^m$ shows that under some conditions the coincidence of all ordinary and symbolic powers of $I$ is equivalent to only having $I^{(m)} = I^m$; see [5, Theorem 2.3]. See also [13, Theorem 2.6] and [29] for some various results on how the equality of finite many symbolic and ordinary powers of an ideal may imply their equality for all powers. Turning back to binomial ideals in the combinatorial context, in [10, Theorem 4.1] and [24, Theorem 3.11], for the binomial edge ideal $J$ of the graph $G$ under some conditions it is shown that $J^{(2)} = J^2$ yields $J^{(m)} = J^m$ for all $m$. It is asked in [10, Question 5.6] whether the same implication holds for every simple graph $G$. On the other hand, if $G$ is bipartite, then the binomial edge ideal of $G$ is isomorphic to $\mathcal{I}_G$ as clarified in [27, Remark 3.1]. Furthermore, suppose that one of the following conditions holds:

- $G$ has a net as an induced subgraph;
- Clique number of $G$ is greater than 3;
- $G$ has a PT as an induced subgraph.

While we will see in the next section that the symbolic and ordinary powers of $\mathcal{I}_G$ for these graphs do not coincide for all powers, the argument used in the proof of these results shows that these powers diverge at the second power. So we raise to ask the following question:

**Question 1.1.** (Question 3.1) *Let $\mathcal{I}_G$ be the parity binomial edge ideal of a graph $G$. Is it true that the following conditions are equivalent:*

*(i) $\mathcal{I}_G^{(2)} = \mathcal{I}_G^2$;*

*(ii) $\mathcal{I}_G^{(t)} = \mathcal{I}_G^t$ for every $t \geq 1$.*

See also [12, Problem 5.14]. Finally our observations through computer experiments, and arguing for some graphs with a small number of vertices persuade us to conjecture a generalization of Theorem 2.6 for so-called DO-free graphs as stated in Conjecture 3.2.

The authors would like to thank the referee for helpful suggestions and inspiring remarks.

## 2 Symbolic Powers of Parity Binomial Edge Ideals

Throughout the paper, $G$ is a simple graph on the vertex set $V(G) = [n]$ with the edge set $E(G)$, and $S = K[x_1, \ldots, x_n, y_1, \ldots, y_n]$ is the polynomial ring over a field $K$ with char $K \neq 2$ in $2n$ variables. For each monomial $u = \mathbf{x}^{\mathbf{a}}\mathbf{y}^{\mathbf{b}} = x_1^{a_1}\ldots x_n^{a_n} y_1^{b_1}\ldots y_n^{b_n} \in S$, we set

$$\deg_{\mathbf{x}} u = \sum_{i=1}^{n} a_i, \tag{1}$$

and

$$\deg_{\mathbf{y}} u = \sum_{i=1}^{n} b_i. \tag{2}$$

It is usual to set $f_{ij} = x_i y_j - x_j y_i$ and $g_{ij} = x_i x_j - y_i y_j$. The parity binomial edge ideal of $G$ denoted by $\mathcal{I}_G$ is generated by those binomials $g_{ij}$ with $\{i,j\} \in E(G)$. Moreover, the number of connected components of $G$ is denoted by $c(G)$, the number of bipartite connected components by $c_0(G)$, and the number of non-bipartite connected components by $c_1(G)$. An $(i,j)$-walk of length $r-1$ in $G$ is a sequence $i = i_1, i_2, \ldots, i_r = j$ of its vertices such that $\{i_k, i_{k+1}\} \in E(G)$ for every $k = 1, \ldots, r-1$. The walk is odd (resp. even) if its length is odd (resp. even). Let $\mathcal{S} \subseteq V(G)$. Then $G_{\mathcal{S}}$ denotes the induced subgraph of $G$ on $V(G)\backslash \mathcal{S}$. Consider the saturation

$$\mathcal{J}_G = \mathcal{I}_G : (\prod_{i \in V(G)} x_i y_i)^{\infty}$$

of $\mathcal{I}_G$, and for each non-bipartite connected component $N$ of $G$

$$\mathfrak{p}^{+}(N) = (x_i + y_i : i \in V(N)),$$

and

$$\mathfrak{p}^{-}(N) = (x_i - y_i : i \in V(N)).$$

For the ideal saturation $\mathcal{J}_G$ with the above setting, one has the following result.

**Proposition 2.1.** ([26, Proposition 2.7]) *Let $G$ be a graph. Then*

$$\begin{aligned}\mathcal{J}_G &= (g_{ij} : \text{there is an odd } (i,j)\text{-walk in } G) \\ &+ (f_{ij} : \text{there is an even } (i,j)\text{-walk in } G).\end{aligned}$$

**Proposition 2.2.** ([26, Proposition 4.2]) *Let $G$ be a graph consisting of bipartite connected components $B_1, \ldots, B_{c_0(G)}$ and non-bipartite connected components $N_1, \ldots, N_{c_1(G)}$ and assume that $char(K) \neq 2$, then $\mathcal{J}_G$ is radical, and its minimal primes are the $2^{c_1(G)}$ ideals*

$$\sum_{i=1}^{c_0(G)} \mathcal{J}_{B_i} + \sum_{i=1}^{c_1(G)} \mathfrak{p}^{\sigma_i}(N_i),$$

*where $(\sigma_1, \sigma_2, \ldots, \sigma_{c_1(G)})$ ranges over $\{+,-\}^{c_1(G)}$.*

Let $\mathfrak{s}(G) = c_0(G) + c(G)$. A set $\mathcal{S} \subseteq V(G)$ is called a **disconnector** of $G$ if for every $s \in \mathcal{S}$
$$\mathfrak{s}(G_\mathcal{S}) > \mathfrak{s}(G_{\mathcal{S}\setminus\{s\}}),$$
where $G_\mathcal{S}$ is the subgraph of $G$ on $V(G)\setminus\mathcal{S}$, as introduced before. Let $\mathcal{S} \subseteq V(G)$ be a disconnector of $G$ and for each $s \in \mathcal{S}$, consider the set $\mathcal{C}_{G_\mathcal{S}}(s)$ of those connected components of $G_\mathcal{S}$ which are joined when adding $s$. A minimal prime $\mathfrak{p}$ of $\mathcal{J}_{G_\mathcal{S}}$ is **sign-split** if for each $s \in \mathcal{S}$ which $\mathcal{C}_{G_\mathcal{S}}(s)$ contains no bipartite graphs, the prime summands of $\mathfrak{p}$ corresponding to connected components in $\mathcal{C}_{G_\mathcal{S}}(s)$ are not all equal to $\mathfrak{p}^+$ or all equal to $\mathfrak{p}^-$.

**Theorem 2.3.** ([26, Theorem 4.15]) *The minimal primes of $\mathcal{I}_G$ are the ideals $\mathfrak{m}_\mathcal{S} + \mathfrak{p}$, where $\mathcal{S} \subseteq V(G)$ is a disconnector of $G$ and $\mathfrak{p}$ is a sign-split minimal prime of $\mathcal{J}_{G_\mathcal{S}}$ and*
$$\mathfrak{m}_\mathcal{S} = (x_s, y_s : s \in \mathcal{S}).$$

A graph $G$ is called a **net** if it is of the form of Figure 1. We call a graph net-free if it does not have a net as an induced subgraph.

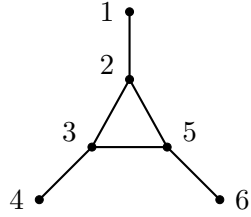

Figure 1

**Theorem 2.4.** *Let $\mathcal{I}_G$ be the parity binomial edge ideal of a graph $G$. If*
$$\mathcal{I}_G^{(t)} = \mathcal{I}_G^t,$$
*for every $t \geq 1$, then $G$ is net-free.*

*Proof.* Assume that $G$ contains a net as an induced subgraph, say $H$ as shown in Figure 1, with edges
$$E(H) = \{\{1,2\}, \{2,3\}, \{3,4\}, \{3,5\}, \{2,5\}, \{5,6\}\}.$$

Set
$$\begin{aligned} h &= y_4 y_6 g_{25} f_{13} - y_3 y_5 g_{16} f_{24} - y_3 y_4 g_{12} f_{56} \\ &= y_1 y_6 g_{23} f_{45} - y_2 y_5 g_{14} f_{36} + y_5 y_6 g_{34} f_{12}. \end{aligned}$$

We will show that $h \in \mathcal{I}_G^{(2)} \setminus \mathcal{I}_G^2$.

We first see that $h \in \mathcal{I}_G^{(2)}$. By [26, Theorem 5.5], the parity binomial edge ideal $\mathcal{I}_G$ is a radical ideal. Hence,
$$\mathcal{I}_G = \bigcap_{\mathfrak{q} \in \mathrm{Min}(\mathcal{I}_G)} \mathfrak{q}.$$

In particular,
$$\mathcal{I}_G^{(2)} = \bigcap_{\mathfrak{q} \in \mathrm{Min}(\mathcal{I}_G)} \mathfrak{q}^{(2)},$$

and regarding the description of minimal prime ideals of $\mathcal{I}_G$ given in Theorem 2.3, it is enough to show that $h \in (\mathfrak{m}_\mathcal{S} + \mathfrak{p})^{(2)}$ for each disconnector $\mathcal{S} \subseteq V(G)$ and each sign-split minimal prime $\mathfrak{p}$ of $\mathcal{J}_{G_\mathcal{S}}$. Thus, suppose that $\mathcal{S} \in V(G)$ is a disconnector of $G$ and consider the following cases:

Case 1. Assume that $\mathcal{S} \cap [6] = \emptyset$. So $G_\mathcal{S}$ has a non-bipartite connected component $N_1$ including $H$ and $\mathfrak{m}_\mathcal{S}$ does not intersect $\{x_i, y_i : i \in [6]\}$. Suppose that $\mathfrak{p}$ is a sign-split minimal prime of $\mathcal{J}_{G_\mathcal{S}}$ and $\mathfrak{p}^{\sigma_1}(N_1)$ is its summand corresponding to the connected component $N_1$; see Proposition 2.2. Consider the following presentations of $f_{ij}$ and $g_{ij}$ for each $i$ and $j$:

$$g_{ij} = x_i(x_j + y_j) - y_j(x_i + y_i) = x_i(x_j - y_j) + y_j(x_i - y_i), \tag{3}$$

$$f_{ij} = y_j(x_i + y_i) - y_i(x_j + y_j) = x_j(x_i - y_i) - x_i(x_j - y_j). \tag{4}$$

By these presentations of $f_{ij}$'s and $g_{ij}$'s, one can see that $h$ has a presentation in which each summand is divisible by two binomials of the form $(x_i + y_i)$ (or similarly two binomials of the form $(x_i - y_i)$ using the second presentation of $f_{ij}$'s and $g_{ij}$'s given above) with $i \in V(H) \subseteq V(N_1)$. Hence,

$$h \in (\mathfrak{p}^{\sigma_1}(N_1))^2 \subseteq \mathfrak{p}^2 \subseteq (\mathfrak{m}_\mathcal{S} + \mathfrak{p})^2.$$

Now it is enough to recall that $(\mathfrak{m}_\mathcal{S} + \mathfrak{p})^2 \subseteq (\mathfrak{m}_\mathcal{S} + \mathfrak{p})^{(2)}$.

Case 2. Next assume that $|\mathcal{S} \cap [6]| = 1$, say $\mathcal{S} \cap [6] = \{i\}$.

2.1. Suppose that $G_\mathcal{S}$ has a non-bipartite connected component $N_1$ including some vertices of $H$. Consider the summand $\mathfrak{p}^{\sigma_1}(N_1)$ corresponding to the component $N_1$ of a sign-split minimal prime $\mathfrak{p}$ of $\mathcal{J}_{G_\mathcal{S}}$, as described in Theorem 2.2. If $\mathcal{S} \cap [6] \in \{1, 3, 4, 5, 6\}$, take the first presentation of $h$, namely

$$h = y_4 y_6 g_{25} f_{13} - y_3 y_5 g_{16} f_{24} - y_3 y_4 g_{12} f_{56},$$

and if $\mathcal{S} \cap [6] \in \{2\}$, take its other presentation, that is,

$$h = y_1 y_6 g_{23} f_{45} - y_2 y_5 g_{14} f_{36} + y_5 y_6 g_{34} f_{12}.$$

Now replace $f_{ij}$'s and $g_{ij}$'s in $h$ by their presentation in (3) and (4) whenever $i$ and $j$ are both vertices of $N_1$ (choose the first presentation of $f_{ij}$'s and $g_{ij}$'s for $\sigma_1 = +$ and the second ones for $\sigma_1 = -$). Then one can see that $h$ has a presentation in which each summand is divisible by a product of a binomial $(x_i + y_i)$ (or similarly $(x_i - y_i)$ in the case of $\sigma_1 = -$) with $i \in V(N_1)$ and a variable $x_j$ or $y_j$ with $j \in \mathcal{S} \cap [6]$. In other words, $h$ has a presentation in which each summand is divisible by a product of a binomial $(x_i + y_i)$ (or $(x_i - y_i)$) in $\mathfrak{p}^{\sigma_1}(N_1)$ and a variable $x_j$ or $y_j$ in $\mathfrak{m}_\mathcal{S}$. Hence,

$$\begin{aligned} h &\in \mathfrak{p}^{\sigma_1}(N_1)\mathfrak{m}_\mathcal{S} \subseteq \mathfrak{p}\mathfrak{m}_\mathcal{S} \\ &\subseteq (\mathfrak{m}_\mathcal{S} + \mathfrak{p})^2 \subseteq (\mathfrak{m}_\mathcal{S} + \mathfrak{p})^{(2)}, \end{aligned}$$

as desired.

2.2. If $G_\mathcal{S}$ has no non-bipartite connected component which has common vertices with $H$, then an induced subgraph of $H - i$ is part of a bipartite connected component $B_1$ of $G_\mathcal{S}$. This case only happens if $\mathcal{S} \cap [6] = \{2\}$ or $\{3\}$ or $\{5\}$. If $\mathcal{S} \cap [6] = \{3\}$ or $\{5\}$, take the first presentation of $h$, that is,

$$h = y_4 y_6 g_{25} f_{13} - y_3 y_5 g_{16} f_{24} - y_3 y_4 g_{12} f_{56}.$$

Otherwise, if $\mathcal{S} \cap [6] = \{2\}$, choose the other presentation, namely

$$h = y_1 y_6 g_{23} f_{45} - y_2 y_5 g_{14} f_{36} + y_5 y_6 g_{34} f_{12}.$$

Now one can see that by choosing these presentations of $h$ and setting $i = \mathcal{S} \cap [6]$, in each summand of $h$ either the variable $y_i$ appears or for some $\ell$, $g_{i\ell}$ or $f_{i\ell}$ appear. Hence, $h \in \mathfrak{m}_\mathcal{S}$. On the other hand, by choosing an appropriate presentation of $h$ as stated above, in each summand a binomial $f_{k\ell}$ or $g_{k\ell}$ with $k, \ell \in V(B_1)$ appear such that these binomials are elements of $\mathcal{J}_{B_1}$ by Proposition 2.1. Thus, regarding Proposition 2.2, one obtains that $h$ is also in $\mathfrak{p}$ for each sign-split minimal prime $\mathfrak{p}$ of $\mathcal{J}_{G_\mathcal{S}}$. Finally, regarding these factors appearing in each summand of $h$, one has $h \in \mathfrak{m}_\mathcal{S} \mathfrak{p}$ for each sign-split minimal prime $\mathfrak{p}$ of $\mathcal{J}_{G_\mathcal{S}}$, and this implies that

$$h \in (\mathfrak{m}_\mathcal{S} + \mathfrak{p})^2 \subseteq (\mathfrak{m}_\mathcal{S} + \mathfrak{p})^{(2)}.$$

Case 3. Assume that $|\mathcal{S} \cap [6]| > 1$. For each monomial summand $m$ of $h$ and each $i \in [6]$, we have

$$x_i | m \text{ or } y_i | m.$$

Hence $m \in \mathfrak{m}_\mathcal{S}^2$. Consequently, $h \in \mathfrak{m}_\mathcal{S}^2 \subseteq (\mathfrak{m}_\mathcal{S} + \mathfrak{p})^2$.

Now we are going to show that $h \notin \mathcal{I}_G^2$. Let $\mathcal{I}_H$ denote the parity binomial edge ideal of $H$ in $K[x_1, \ldots, x_6, y_1, \ldots, y_6]$. Applying the same argument as used in the proof of [25, Proposition 3.3], one has $\mathcal{I}_H^t = \mathcal{I}_G^t \cap K[x_1, \ldots, x_6, y_1, \ldots, y_6]$ for every $t \geq 1$. So it is enough to show that $h \notin \mathcal{I}_H^2$. Suppose that $h \in \mathcal{I}_H^2$, and consider its presentation as a finite sum

$$h = \sum_{\mathbf{a}=(a_1, a_2, a_3, a_4)} \mu_\mathbf{a} g_{a_1 a_2} g_{a_3 a_4}, \tag{5}$$

where $\mu_{\mathbf{a}} \in K[x_1, \ldots, x_6, y_1, \ldots, y_6]$, and $\{a_1, a_2\}, \{a_3, a_4\} \in E(H)$. Now, on the one hand, notice that $\deg_{\mathbf{x}} u = 4 - \deg_{\mathbf{y}} u$ can be $0, 2$, or $4$ for each monomial term $u$ of a generator $g_{a_1 a_2} g_{a_3 a_4}$. On the other hand, recall that $x_1 x_4 x_6 y_2 y_3 y_5$ is a monomial term of $h$, and consequently divides a monomial term $v$ of $\mu_{\mathbf{a}} g_{a_1 a_2} g_{a_3 a_4}$ with degree 6 for some $\mathbf{a}$ in the presentation (5). Regarding the fact

$$\deg_{\mathbf{x}} v = 6 - \deg_{\mathbf{y}} v = \deg_{\mathbf{x}} x_1 x_4 x_6 y_2 y_3 y_5 = 3,$$

the set $\{1, 4, 6\}$ must contain the edge $\{a_1, a_2\}$ or $\{a_3, a_4\}$ of the graph $H$. Otherwise, $\deg_{\mathbf{y}} v \geq 4$ which is not the case. Now it is enough to notice that $\{1, 4, 6\}$ does contain no edge of H. By this contradiction, we deduce that $h \notin \mathcal{I}_H^2$. □

A **clique** of a graph $G$ is a complete subgraph of $G$. A **maximal clique** of $G$ is a clique of $G$ that cannot be extended to a larger clique of $G$. The **clique number** of $G$, denoted by $\omega(G)$, is the number of vertices in a maximum clique of $G$.

**Theorem 2.5.** *Let $\mathcal{I}_G$ be the parity binomial edge ideal of a graph $G$. If*

$$\mathcal{I}_G^{(t)} = \mathcal{I}_G^t,$$

*for every $t \geq 1$, then $\omega(G) \leq 3$.*

*Proof.* Assume that $G$ contains $K_4$ with edges

$$E(K_4) = \{\{1, 2\}, \{2, 3\}, \{3, 4\}, \{1, 4\}, \{1, 3\}, \{2, 4\}\}.$$

Set $h = x_1 x_2 x_4 y_3 - x_1 y_2 y_3 y_4 - x_2 x_3 x_4 y_1 + x_3 y_1 y_2 y_4$. Recall that $g_{ij} = x_i x_j - y_i y_j$ and $f_{ij} = x_i y_j - x_j y_i$. So

$$h = g_{24} f_{13}.$$

We will show that $h \in \mathcal{I}_G^{(2)} \setminus \mathcal{I}_G^2$.

We first see that $h \in \mathcal{I}_G^{(2)}$. By the argument used in the proof of Theorem 2.4, it is enough to show that $h \in (\mathfrak{m}_{\mathcal{S}} + \mathfrak{p})^{(2)}$ for each disconnector $\mathcal{S} \subseteq V(G)$ and each sign-split minimal prime $\mathfrak{p}$ of $\mathcal{J}_{G_{\mathcal{S}}}$. So suppose that $\mathcal{S} \in V(G)$ is a disconnector of $G$, and consider the following cases:

Case 1. Assume that $\mathcal{S} \cap [4] = \emptyset$. So $G_{\mathcal{S}}$ has a non-bipartite connected component $N_1$ including $K_4$ and $\mathfrak{m}_{\mathcal{S}}$ does not include $\{x_i, y_i : i \in [4]\}$. By Theorem 2.3, we know that each minimal prime of $\mathcal{I}_G$ has a summand $\mathfrak{p}$ which is a sign-split minimal prime of $\mathcal{J}_{G_{\mathcal{S}}}$ for some disconnector $\mathcal{S} \subseteq V(G)$ of $G$.

Considering $\mathfrak{p}$ as introduced above, we discuss the summand $\mathfrak{p}^{\sigma_1}$ of $\mathfrak{p}$. By the following presentations of $h$ given below, we will see that each summand of $h$ is of the form of a product of two binomials of the form $(x_i + y_i)$ or two binomials of the form $(x_i - y_i)$. So we conclude that $h \in \mathfrak{p}^2 \subseteq (\mathfrak{m}_{\mathcal{S}} + \mathfrak{p})^2$.

1.1. If $\sigma_1 = +$, applying
$$g_{24} = x_2(x_4 + y_4) - y_4(x_2 + y_2),$$
and
$$f_{13} = y_3(x_1 + y_1) - y_1(x_3 + y_3),$$
one obtains the presentation
$$\begin{aligned} h &= x_2 y_3(x_1 + y_1)(x_4 + y_4) - x_2 y_1(x_3 + y_3)(x_4 + y_4) - \\ & \quad y_3 y_4(x_1 + y_1)(x_2 + y_2) + y_1 y_4(x_2 + y_2)(x_3 + y_3). \end{aligned}$$

1.2. If $\sigma_1 = -$, applying
$$g_{24} = x_2(x_4 - y_4) + y_4(x_2 - y_2),$$
and
$$f_{13} = x_3(x_1 - y_1) - x_1(x_3 - y_3),$$
one obtains the presentation
$$\begin{aligned} h &= x_2 x_3(x_1 - y_1)(x_4 - y_4) - x_1 x_2(x_3 - y_3)(x_4 - y_4) + \\ & \quad x_3 y_4(x_1 - y_1)(x_2 - y_2) - x_1 y_4(x_2 - y_2)(x_3 - y_3). \end{aligned}$$

Case 2. Assuming that $|\mathcal{S} \cap [4]| = 1$, we know that $G_\mathcal{S}$ has a non-bipartite connected component $N_1$ including some vertices of $K_4$. Now in the following presentations of $h$, each summand is divisible by a product of a binomial of the form $(x_i + y_i)$ or $(x_i - y_i)$ with $i \in V(N_1)$ and a variable $x_j$ or $y_j$ with $j$ in the disconnector $\mathcal{S}$. This implies that $h \in (\mathfrak{m}_\mathcal{S} + \mathfrak{p})^2$.

2.1. If $\sigma_1 = +$, and $\mathcal{S} \cap [4] = \{1\}$ or $\{3\}$, then
$$\begin{aligned} h &= f_{13}[x_2(x_4 + y_4) - y_4(x_2 + y_2)] \\ &= x_1 y_3[x_2(x_4 + y_4) - y_4(x_2 + y_2)] - x_3 y_1[x_2(x_4 + y_4) - y_4(x_2 + y_2)]. \end{aligned}$$

2.2. If $\sigma_1 = +$, and $\mathcal{S} \cap [4] = \{2\}$ or $\{4\}$, then
$$\begin{aligned} h &= g_{24}[y_3(x_1 + y_1) - y_1(x_3 + y_3)] \\ &= x_2 x_4[y_3(x_1 + y_1) - y_1(x_3 + y_3)] - y_2 y_4[y_3(x_1 + y_1) - y_1(x_3 + y_3)]. \end{aligned}$$

2.3. If $\sigma_1 = -$, and $\mathcal{S} \cap [4] = \{1\}$ or $\{3\}$, then
$$\begin{aligned} h &= f_{13}[x_2(x_4 - y_4) + y_4(x_2 - y_2)] \\ &= x_1 y_3[x_2(x_4 - y_4) + y_4(x_2 - y_2)] - x_3 y_1[x_2(x_4 - y_4) + y_4(x_2 - y_2)]. \end{aligned}$$

2.4. If $\sigma_1 = -$, and $\mathcal{S} \cap [4] = \{2\}$ or $\{4\}$, then
$$\begin{aligned} h &= g_{24}[x_3(x_1 - y_1) - x_1(x_3 - y_3)] \\ &= x_2 x_4[x_3(x_1 - y_1) - x_1(x_3 - y_3)] - y_2 y_4[x_3(x_1 - y_1) - x_1(x_3 - y_3)]. \end{aligned}$$

Case 3. Finally, assume that $|\mathcal{S} \cap [4]| > 1$. For each monomial summand $m$ of $h$, we have $x_i|m$ or $y_i|m$ for every $i \in [4]$. So $m \in \mathfrak{m}_\mathcal{S}^2$, and consequently $h \in (\mathfrak{m}_\mathcal{S} + \mathfrak{p})^2$.

As we have seen, $h \in (\mathfrak{m}_\mathcal{S} + \mathfrak{p})^2 \subseteq (\mathfrak{m}_\mathcal{S} + \mathfrak{p})^{(2)}$ in all cases above. Thus, $h \in \mathcal{I}_G^{(2)}$.

Now we show that $h \notin \mathcal{I}_G^2$. Since $\mathcal{I}_G$ is induced by binomials of type $x_i x_j - y_i y_j$, we conclude that $\deg_\mathbf{x} u$ and $\deg_\mathbf{y} u$, as defined in (1) and (2), are even numbers for each monomial term $u$ of a polynomial in the minimal set of homogeneous generators $\{g_{ij}g_{k\ell} : \{i,j\}, \{k,\ell\} \in E(G)\}$ of $\mathcal{I}_G^2$. Hence, for the homogeneous polynomial $h$ of degree 4, we deduce that $h \notin \mathcal{I}_G^2$. $\square$

We call a graph a **PT** if it is a pentagon and a triangle graph sharing an edge; see Figure 2. A graph is **PT-free** if it does not contain a PT as an induced subgraph.

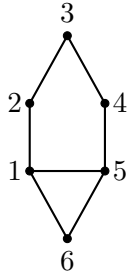

Figure 2

**Theorem 2.6.** *Let $\mathcal{I}_G$ be the parity binomial edge ideal of a graph $G$. If*
$$\mathcal{I}_G^{(t)} = \mathcal{I}_G^t,$$
*for every $t \geq 1$, then $G$ is PT-free.*

*Proof.* Assume that $G$ has an induced PT subgraph, say H on 6 vertices as shown in Figure 2, with edges
$$E(H) = \{\{1,2\}, \{2,3\}, \{3,4\}, \{4,5\}, \{5,6\}, \{1,6\}\} \cup \{\{1,5\}\}.$$

Set the polynomial $h$ as follows:
$$h = y_1 y_2 g_{36} f_{45}.$$

We will show that $h \in \mathcal{I}_G^{(2)} \setminus \mathcal{I}_G^2$.

We first see that $h \in \mathcal{I}_G^{(2)}$. Recall that by the argument used in the proof of Theorem 2.4, it is enough to show that $h \in (\mathfrak{m}_\mathcal{S} + \mathfrak{p})^{(2)}$ for each disconnector $\mathcal{S} \subseteq V(G)$ and each sign-split minimal prime $\mathfrak{p}$ of $\mathcal{J}_{G_\mathcal{S}}$. So suppose that $\mathcal{S} \in V(G)$ is a disconnector of $G$, and consider the following cases:

Case 1. Assume that $\mathcal{S} \cap [6] = \emptyset$. So $G_\mathcal{S}$ has a non-bipartite connected component $N_1$ including $H$ and $\mathfrak{m}_\mathcal{S}$ does not intersect $\{x_i, y_i : i \in [6]\}$. Suppose that $\mathfrak{p}$ is a sign-split minimal prime of $\mathcal{J}_{G_\mathcal{S}}$ and $\mathfrak{p}^{\sigma_1}(N_1)$ is its summand corresponding to the connected component $N_1$; see Proposition 2.2. Consider the following presentations of $f_{45}$ and $g_{36}$:

$$g_{36} = x_3(x_6 + y_6) - y_6(x_3 + y_3) = x_3(x_6 - y_6) + y_6(x_3 - y_3), \tag{6}$$

$$f_{45} = y_5(x_4 + y_4) - y_4(x_5 + y_5) = x_5(x_4 - y_4) - x_4(x_5 - y_5). \tag{7}$$

By these presentations of $f_{45}$ and $g_{36}$, one can see that $h$ has a presentation in which each summand is divisible by two binomials of the form $(x_i + y_i)$ (or similarly two binomials of the form $(x_i - y_i)$ using the second presentation of $f_{45}$ and $g_{36}$ given above) with $i \in V(H) \subseteq V(N_1)$. Hence,

$$h \in (\mathfrak{p}^{\sigma_1}(N_1))^2 \subseteq \mathfrak{p}^2 \subseteq (\mathfrak{m}_\mathcal{S} + \mathfrak{p})^2.$$

Now it is enough to recall that $(\mathfrak{m}_\mathcal{S} + \mathfrak{p})^2 \subseteq (\mathfrak{m}_\mathcal{S} + \mathfrak{p})^{(2)}$.

Case 2. Next assume that $|\mathcal{S} \cap [6]| = 1$, say $\mathcal{S} \cap [6] = \{i\}$.

2.1. First suppose that $G_\mathcal{S}$ has a non-bipartite connected component $N_1$ including vertices $V(H) \setminus \{i\}$. Using the presentation of $f_{45}$ and $g_{36}$ given in (6) and (7), one obtains a presentation of $h$ in which each summand is divisible by a product of the variable $x_i$ or $y_i$ and a binomial $x_k + y_k$ (or similarly $x_k - y_k$) with $k \in V(H) \setminus \{i\} \subseteq V(N_1)$. Hence, $h \in (\mathfrak{m}_\mathcal{S})(\mathfrak{p}^+(N_1)) \subseteq (\mathfrak{m}_\mathcal{S} + \mathfrak{p}^+(N_1))^2$ (or similarly $h \in (\mathfrak{m}_\mathcal{S} + \mathfrak{p}^-(N_1))^2$). This guarantees that for each sign-split minimal prime $\mathfrak{p}$ of $\mathcal{J}_{G_\mathcal{S}}$, we have $h \in (\mathfrak{m}_\mathcal{S} + \mathfrak{p})^2 \subseteq (\mathfrak{m}_\mathcal{S} + \mathfrak{p})^{(2)}$.

2.2. If $G_\mathcal{S}$ has no non-bipartite connected component that has common vertices with $H$, then $H - i$ is in a bipartite connected component $B_1$ of $G_\mathcal{S}$. This case only happens if $\mathcal{S} \cap [6] = \{1\}$ or $\{5\}$. Without loss of generality, assume that $\mathcal{S} \cap [6] = \{1\}$. By Proposition 2.1 we have

$$g_{36} \in \mathcal{J}_{B_1}.$$

This yields that

$$f_{45} g_{36} \in \mathfrak{m}_\mathcal{S} \mathcal{J}_{B_1}.$$

Hence, $h \in \mathfrak{m}_\mathcal{S} \mathcal{J}_{B_1}$. By Proposition 2.2, for each sign-split minimal prime $\mathfrak{p}$ of $\mathcal{J}_{G_\mathcal{S}}$, we obtain that $h \in \mathfrak{m}_\mathcal{S} \mathfrak{p}$. So $h \in (\mathfrak{m}_\mathcal{S} + \mathfrak{p})^2 \subseteq (\mathfrak{m}_\mathcal{S} + \mathfrak{p})^{(2)}$.

Case 3. Finally, assume that $|\mathcal{S} \cap [6]| > 1$. Then each term $m$ of $h$ is divisible by $x_i$ or $y_i$ for every $i \in [6]$. So $m \in \mathfrak{m}_\mathcal{S}^2$, and $h \in (\mathfrak{m}_\mathcal{S} + \mathfrak{p})^2$.

Now we are going to show that $h \notin \mathcal{I}_G^2$. Let $\mathcal{I}_H$ denote the parity binomial edge ideal of $H$ in $K[x_1, \ldots, x_6, y_1, \ldots, y_6]$. Applying the same argument as used in the proof of [25, Proposition 3.3], one has $\mathcal{I}_H^t = \mathcal{I}_G^t \cap K[x_1, \ldots, x_6, y_1, \ldots, y_6]$ for every $t \geq 1$. So it is enough to show that $h \notin \mathcal{I}_H^2$ which can be computed, for example, by Macualay2 [16]. $\square$

# 3  Open Problems

Let $G$ be a simple graph. Suppose that $J$ is the binomial edge ideal of $G$, a binomial ideal that is isomorphic to $\mathcal{I}_G$ when $G$ is a bipartite graph as pointed out in [27, Remark 3.1]; see [19] for the definition and more explanation on binomial edge ideals. It is shown in [10, Theorem 4.1] and [24, Theorem 3.11] that if $G$ is a block graph with Cohen-Macaulay binomial edge ideal or a generalized caterpillar graph for which $J^{(m)} \neq J^m$ for some $m$, then $J^{(2)} \neq J^2$. Ene et al. asked in [10, Question 5.6] whether $J^{(2)} = J^2$ implies that $J^{(m)} = J^m$ for all $m$. A positive answer to this question implies the same consequence for the parity binomial edge ideal of bipartite graphs. On the other hand, we have seen in the previous section that the symbolic and ordinary power of the parity binomial edge ideal diverge at second power when $G$ has a net or PT as an induced subgraph, or clique number of $G$ is greater than 3. So we ask the following question:

**Question 3.1.** *Let $\mathcal{I}_G$ be the parity binomial edge ideal of a graph $G$. Is it true that the following conditions are equivalent:*

(i) $\mathcal{I}_G^{(2)} = \mathcal{I}_G^2$;

(ii) $\mathcal{I}_G^{(t)} = \mathcal{I}_G^t$ for every $t \geq 1$.

A combinatorial description of minimal prime ideals of $\mathcal{I}_G$ is presented by Kahle et al. (see Theorem 2.3). It could be of interest to find a combinatorial description based on underlying graph $G$ for embedded prime ideals of higher powers of $\mathcal{I}_G$. Besides, while for the graphs we consider in the previous section, the divergence of symbolic and ordinary powers of $\mathcal{I}_G$ occurs at second power, it is unclear in which higher powers this divergence would repeat.

By a result of Brodmann in [4], the set $\mathrm{Ass}(R/I^m)$ stabilizes for large enough powers of each ideal $I$ in a Noetherian ring $R$. The smallest number $k$ for which $\mathrm{Ass}(R/I^m) = \mathrm{Ass}(R/I^k)$ for every $m \geq k$ is called the index of stability of $I$. Let $k$ be the index of stability of $I$. As Morey pointed out in [29], the result of Brodmann guarantees that if $\mathrm{Ass}(R/I^k) = \mathrm{Min}(I)$, then eventually $I^{(m)} = I^m$. It is also clear that eventually $I^{(m)} \neq I^m$ if $\mathrm{Ass}(R/I^k) \neq \mathrm{Min}(I)$. What is the least number $\ell$ such that the comparison of $I^{(\ell)}$ and $I^\ell$ can determine the eventual behavior of these powers? Is there a combinatorial description for such a number?

We call a graph **double-odd cycle**, or **DO** for short, if it is two odd cycles sharing an edge, See Figure 3. A graph is **DO-free** if it does not contain a DO as an induced subgraph. In particular, PT-free graphs are among this class of graphs. Using Macualay2 [16], computer experiments for a double-odd cycle $G$ with a small number of vertices show that the second symbolic and ordinary power of $\mathcal{I}_G$ are not equal.

**Conjecture 3.2.** *Let $\mathcal{I}_G$ be the parity binomial edge ideal of a graph $G$. If we have*

$$\mathcal{I}_G^{(t)} = \mathcal{I}_G^t$$

*for every $t \geq 1$, then $G$ is DO-free.*

Figure 3

In support of our conjecture, assume that $G$ has the induced DO subgraph H on $4k$ vertices where $H$ is a cycle of length $4k$ including the diagonal $\{1, 3\}$, say

$$E(H) = \{\{1,2\}, \{2,3\}, \ldots, \{4k-1, 4k\}, \{4k, 1\}\} \cup \{\{1,3\}\}.$$

Set the polynomial $h$ as follows:

$$h = y_4 y_5 \ldots \widehat{y_{2k+2}} y_{2k+3} \ldots y_{4k} g_{13} f_{2\,2k+2}.$$

Regarding Conjecture 3.2, the polynomial $h$ is the candidate to be in $\mathcal{I}_G^{(2)} \setminus \mathcal{I}_G^2$. Using Macualay2 [16], computation shows that $h \notin \mathcal{I}_G^2$ for small numbers of $k$. On the other hand, one can see that $h \in \mathcal{I}_G^{(2)}$ as follows:

Since
$$\mathcal{I}_G^{(2)} = \bigcap_{\mathfrak{q} \in \mathrm{Min}(\mathcal{I}_G)} \mathfrak{q}^{(2)},$$

using the description of minimal primes of $\mathcal{I}_G$ given in Theorem 2.3, we only need to show that $h \in (\mathfrak{m}_\mathcal{S} + \mathfrak{p})^{(2)}$ for each disconnector $\mathcal{S} \subseteq V(G)$ and each sign-split minimal prime $\mathfrak{p}$ of $\mathcal{J}_{G_\mathcal{S}}$. Let $\mathcal{S} \in V(G)$ be a disconnector of $G$ and consider the following cases:

Case 1. Assume that $\mathcal{S} \cap [4k] = \emptyset$. Then $G_\mathcal{S}$ has a non-bipartite connected component $N_1$ including $H$. Suppose that $\mathfrak{p}$ is a sign-split minimal prime of $\mathcal{J}_{G_\mathcal{S}}$ and $\mathfrak{p}^{\sigma_1}(N_1)$ is its summand corresponding to the connected component $N_1$; see Proposition 2.2. Consider the following presentations of $g_{13}$ and $f_{2\,2k+2}$:

$$g_{13} = x_1(x_3 + y_3) - y_3(x_1 + y_1) = x_1(x_3 - y_3) + y_3(x_1 - y_1), \tag{8}$$

$$f_{2\,2k+2} = y_{2k+2}(x_2 + y_2) - y_2(x_{2k+2} + y_{2k+2}) = x_{2k+2}(x_2 - y_2) - x_2(x_{2k+2} - y_{2k+2}). \tag{9}$$

Then $h$ has a presentation in which each summand is divisible by a product of two binomials of the form $(x_i + y_i)$ by first presentations in (8) and (9) with $i \in V(H) \subseteq V(N_1)$, and similarly using the second presentations, there exists a presentation of $h$ whose summands are divisible by product of two binomials of the form $(x_i - y_i)$ with $i \in V(H) \subseteq V(N_1)$. Hence

$$h \in (\mathfrak{p}^{\sigma_1}(N_1))^2 \subseteq \mathfrak{p}^2 \subseteq (\mathfrak{m}_\mathcal{S} + \mathfrak{p})^2$$

where $\sigma \in \{-, +\}$. Now recall that $(\mathfrak{m}_\mathcal{S} + \mathfrak{p})^2 \subseteq (\mathfrak{m}_\mathcal{S} + \mathfrak{p})^{(2)}$.

Case 2. Next assume that $|\mathcal{S} \cap [4k]| = 1$, say $\mathcal{S} \cap [4k] = \{i\}$.

2.1. First suppose that $G_\mathcal{S}$ has a non-bipartite connected component $N_1$ including vertices $V(H) \setminus \{i\}$. Using the presentation of $g_{13}$ and $f_{2\,2k+2}$ given in (8) and (9), one obtains a presentation of $h$ in which each summand is divisible by a product of the variable $x_i$ or $y_i$ and a binomial $x_j + y_j$ (or similarly $x_j - y_j$) with $j \in V(H) \setminus \{i\} \subseteq V(N_1)$. Hence $h \in (\mathfrak{m}_\mathcal{S})(\mathfrak{p}^+(N_1)) \subseteq (\mathfrak{m}_\mathcal{S} + \mathfrak{p}^+(N_1))^2$ (or similarly $h \in (\mathfrak{m}_\mathcal{S} + \mathfrak{p}^-(N_1))^2$). As a result, $h \in (\mathfrak{m}_\mathcal{S} + \mathfrak{p})^2$ for each sign-split minimal prime $\mathfrak{p}$ of $\mathcal{J}_{G_\mathcal{S}}$.

2.2. If $G_\mathcal{S}$ has no non-bipartite connected component which has common vertices with $H$, then $H - i$ is in a bipartite connected component $B_1$ of $G_\mathcal{S}$. This happens only if $\mathcal{S} \cap [4k] = \{1\}$ or $\{3\}$. Without loss of generality assume that $\mathcal{S} \cap [4k] = \{1\}$. Considering the even walk $2, 3, \ldots, 2k+1, 2k+2$ of length $2k$, by Proposition 2.1 we have

$$f_{2\,2k+2} \in \mathcal{J}_{B_1}.$$

This yields that

$$g_{13} f_{2\,2k+2} \in \mathfrak{m}_\mathcal{S} \mathcal{J}_{B_1}.$$

Hence $h \in \mathfrak{m}_\mathcal{S} \mathcal{J}_{B_1}$. In particular by Proposition 2.2, $h \in \mathfrak{m}_\mathcal{S} \mathfrak{p}$ for each sign-split minimal prime $\mathfrak{p}$ of $\mathcal{J}_{G_\mathcal{S}}$. Therefore, $h \in (\mathfrak{m}_\mathcal{S} + \mathfrak{p})^2$.

Case 3. Finally assume that $|\mathcal{S} \cap [4k]| > 1$. Then each term $m$ of $h$ is divisible by $x_i$ or $y_i$ for every $i \in [4k]$. So $m \in \mathfrak{m}_\mathcal{S}^2$, and $h \in (\mathfrak{m}_\mathcal{S} + \mathfrak{p})^2$.

## References

[1] A. Alilooee, A. Banerjee. *Packing properties of cubic square-free monomial ideals.* Journal of Algebraic Combinatorics, 54(3), 803-813. 2021.

[2] C. Bahiano. *Symbolic powers of edge ideals.* Journal of Algebra, 273, pp. 517-537, 2004.

[3] S. Bayati, F. Rahmati. *Squarefree vertex cover algebras.* Communications in Algebra, 42(4), pp. 1518-1538, 2014.

[4] M. Brodmann. *Asymptotic stability of* $\mathrm{Ass}(M/I^n M)$. Proceedings of the American Mathematical Society. 74, 16-18, 1979.

[5] S. Cooper, G. Fatabbi, E. Guardo, A. Lorenzini, J. Migliore, U. Nagel, A. Seceleanu, J. Szpond, A. Van Tuyl. *Symbolic powers of codimension two Cohen-Macaulay ideals.* Communications in Algebra, 48(11), 2020.

[6] H. Dao, A. De Stefani, E. Grifo, C. Huneke, L. Núñez-Betancourt. *Symbolic powers of ideals.* In Singularities and foliations. geometry, topology and applications, volume 222 of Springer Proc. Math. Stat., 387-432. Springer, Cham, 2018.

[7] L. Dupont, R. Villarreal. *Edge ideals of clique clutters of comparability graphs and the normality of monomial ideals.* Mathematica Scandinavica, 106(1), 88-98, 2010.

[8] L. Dupont, R. Villarreal. *Symbolic Rees algebras, vertex covers and irreducible representations of Rees cones.* Algebra and Discrete Mathematics, 10(2), 64-86, 2010.

[9] D. Eisenbud. *Commutative Algebra, with a view towards algebraic geometry.* Graduate Texts in Math. 150, Springer-Verlag, New York, 1995.

[10] V. Ene, G. Rinaldo, N. Terai. *Powers of binomial edge ideals with quadratic gröbner bases.* Nagoya Mathematical Journal, pp. 1-23, 2021.

[11] V. Ene, J. Herzog. *On the symbolic powers of binomial edge ideals.* Combinatorial Structures in Algebra and Geometry, pp. 43-50, 2018.

[12] C. Francisco, H. T. Hà, J. Mermin. *Powers of square-free monomial ideals and combinatorics.* Commutative algebra, 373-392, Springer, New York, 2013.

[13] F. Galetto, A. V. Geramita, Y. Shin, A. Van Tuyl. *The symbolic defect of an ideal.* Journal of Pure and Applied Algebra, 223(6), 2709-2731, 2019.

[14] I. Gitler, E. Reyes, R. Villarreal. *Blowup algebras of square-free monomial ideals and some links to combinatorial optimization problems.* Rocky Mountain Journal of Mathematics, 39(1), 71-102, 2009.

[15] I. Gitler, C. Valencia, R. Villarreal. *A note on Rees algebras and the MFMC property.* Beiträge zur Algebra und Geometrie, 48(1), 141-150, 2007.

[16] D. Grayson, M. Stillman. *Macaulay2, a software system for research in algebraic geometry.* http://www.math.uiuc.edu/Macaulay2/.

[17] E. Grifo. *Symbolic Powers.* Lecture Notes for Escuela de Outonõ en Álgebra Conmutativa, 2019. Available at: eloisagrifo.github.io/SymbolicPowersCIMAT.pdf

[18] H. T. Hà, S. Morey. *Embedded associated primes of powers of square-free monomial ideals.* Journal of Pure and Applied Algebra, 214(4), 301-308, 2010.

[19] J. Herzog, T. Hibi, F. Hreinsdóttir, T. Kahle, J. Rauh. *Binomial edge ideals and conditional independence statements.* Advances in Applied Mathematics, 45(3), pp. 317-333, 2010.

[20] J. Herzog, T. Hibi, N. V. Trung. *Symbolic powers of monomial ideals and vertex cover algebras.* Advances in Mathematics, 210(1), pp. 304-322, 2007.

[21] J. Herzog, T. Hibi, N. Trung, X. Zheng. *Standard graded vertex cover algebras, cycles and leaves.* Transactions of the American Mathematical Society, 360(12), pp. 6231-6249, 2008.

[22] J. Herzog, A. Macchia, S. Saeedi Madani, V. Welker. *On the ideal of orthogonal representations of a graph in $\mathbb{R}^2$.* Advances in Applied Mathematics, 71, pp. 146-173, 2015.

[23] C. Huneke. *The primary components of and integral closures of ideals in 3-dimensional regular local rings.* Mathematische Annalen, 275(4), 617-635, 1986.

[24] I. Jahani, S. Bayati, F. Rahmati. *On the Equality of Symbolic and Ordinary Powers of Binomial Edge Ideals.* Bulletin of the Brazilian Mathematical Society, New Series 54(7), 2023.

[25] A. Jayanthan, A. Kumar, R. Sarkar. *Regularity of powers of quadratic sequences with applications to binomial ideals.* Journal of Algebra, 564, pp. 98-118, 2020.

[26] T. Kahle, C. Sarmiento, T. Windisch. *Parity binomial edge ideals.* Journal of Algebraic Combinatorics, 44(1), pp. 99-117, 2016.

[27] A. Kumar. *Lovász–Saks–Schrijver ideals and parity binomial edge ideals of graphs.* European Journal of Combinatorics, 93, 103274, 19, 2021.

[28] J. Montaño, L. Núñez-Betancourt. *Splittings and symbolic powers of square-free monomial ideals.* International Mathematics Research Notices, no. 3, 2304-2320, 2021.

[29] S. Morey. *Stability of associated primes and equality of ordinary and symbolic powers of ideals.* Communications in Algebra, 27(7), 3221-3231, 1999.

[30] M. Ohtani. *Binomial edge ideals of complete multipartite graphs.* Communications in Algebra, 41(10), pp. 3858-3867, 2013.

[31] M. Ohtani. *Graphs and ideals generated by some 2-minors.* Communications in Algebra, 39(3), pp. 905-917, 2011.

[32] A. Simis, W. Vasconcelos, R. Villarreal. *On the ideal theory of graphs.* Journal of Algebra, 167, 389-416, 1994.

[33] J. Sidman, S. Sullivant. *Prolongations and computational algebra.* Canadian Journal of Mathematics, 61, 930-949, 2009.

[34] A. Li, I. Swanson. *Symbolic powers of radical ideals.* The Rocky Mountain Journal of Mathematics, 36, 997-1009, 2006.

[35] T. Szemberg, J. Szpond. *On the containment problem.* Rendiconti del Circolo Matematico di Palermo Series 2, 66(2), 233-245, 2017.

[36] N. V. Trung, T. M. Tuan. *Equality of ordinary and symbolic powers of Stanley-Reisner ideals.* Journal of Algebra, 328(1), pp. 77-93, 2011.

Department of Mathematics and Computer Science, Amirkabir University of Technology (Tehran Polytechnic), Iran

*E-mail address*: `ntaghipour@aut.ac.ir`
*E-mail address*: `shamilabayati@gmail.com`
*E-mail address*: `frahmati@aut.ac.ir`